\documentclass{article}
\usepackage{amsfonts}
\usepackage{multicol}
\author{Jana Archibald}
\title{The Weight System of the Multivariable Alexander Polynomial}
\date{\today}
\newtheorem{theorem}{Theorem}
\newtheorem{lemma}{Lemma}

\newtheorem{definition}{Definition}

\usepackage{graphicx} 
\usepackage{amsmath}

\begin{document}
\maketitle
\begin{abstract}

We derive a formula for the weight system of the multivariable Alexander polynomial using determinants, show that it obeys known relations, and satisfies some of the same relations as the single variable polynomial.

\end{abstract}

\section{Introduction}

When Alexander first defined  the multivariable Alexander polynomial (MVA) \cite{alex} it was defined only up to sign and powers of the variables.  It was shown by Murakami \cite{Mur} that once a normalization and an appropriate change of variables is chosen, the polynomial is of finite type.  In \cite{Mur} a recursive definition of this weight system was given, whose computation would take exponential time.  We give a closed form formula using determinants, which would be computable in polynomial time.

This explicit formula makes it possible to verify that the relations from  \cite{gl}  which hold for the single variable Alexander polynomial also hold for the multivariable Alexander polynomial.

\section{Definition of the MVA}

 \begin{definition} A coloured link (or chord diagram) is a link along with  a   variable associated to each of the components.
\end{definition}

We  will use the normalization of the multivariable Alexander polynomial   given in \cite{Mur}.  This  uses the matrix given by Fox free calculus on the Wirtinger presentation of the link, with the crossings labeled by the arc starting at that crossing.
To compute the multivariable Alexander polynomial one first labels all of the arcs, and then labels the crossings by the exiting lower strand.  This allows us to create a matrix with rows and colums indexed by the arcs.  The $c^{th}$ row is given below.

\begin{multicols}{2}

\begin{center}

\includegraphics[width=35mm]{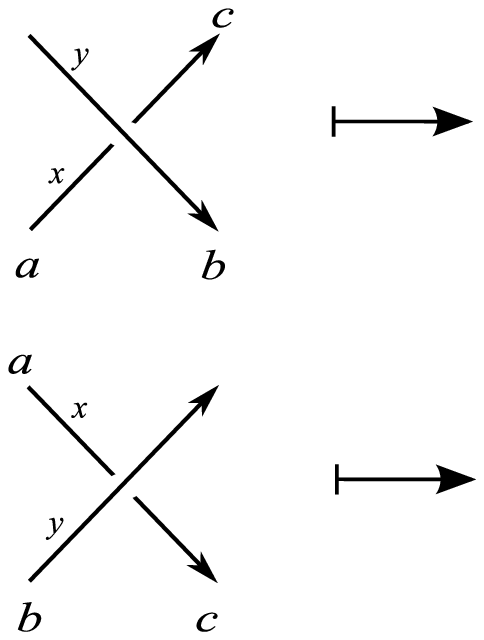}  

\end{center}

 \

 \vspace{2mm}
 
\begin{tabular}{|c|c c c| } \hline
& $a$ & $b$ & $c$  \\  \hline
$c$ & $-1$ & $1-x$ & $y$ \\  \hline
\end{tabular}

\vspace{1.5cm}

\begin{tabular}{|c|c c c| } \hline
& $a$ & $b$ & $c$  \\  \hline
$c$ & $-y$ & $x-1$ & $1$ \\  \hline
\end{tabular}

\end{multicols}

Let $M_i^j$ denote the matrix obtained form $M$ by deleting   the $i^{th}$ row and $j^{th}$ column.  Let $rot(k)$ denote the rotation number of the $k^{th}$ component of the link.  $\mu (k)$ be the number of times that the $k^{th}$ component is the over strand in a crossing.  Let $w_i$ be a word corresponding to a  path from a point to the right of the  $a_i^{th}$ crossing to the unbounded region of the plane.

 \begin{definition} The Multivariable Alexander Polynomial can be defined  as follows,

$$\bigtriangleup (L) = \frac{(-1)^{i+j} \det(M_i^j)}{w_i(t_i-1)}\prod _k t_k^{\frac{rot(k)-\mu (k)}{2}}.$$

\end{definition} 
Note that for links $t_i-1$ divides $M_i^j$ so this is indeed a polynomial, and  that for knots this differs from the usual definition by a factor of $t-1$.   

\

\noindent{\bf Example}

For the following labeled link, we construct the given matrix.   

\begin{multicols}{2}

\includegraphics[width=50mm]{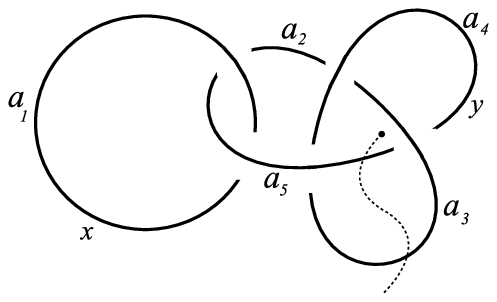}

$\begin{array}{|c | c c c c c|}  \hline
 & a_1 & a_2 & a_3 & a_4 & a_5  \\  \hline
 a_1 & y-1 & 0 & 0 & 0 & 1-x\\
 a_2 & 1-y & x & 0 & 0 & -1 \\
 a_3 & -y & 0& 1 & y-1 & 0 \\
 a_4 & 0 & 0 & -y & 1 & y-1 \\
 a_5 & 0 & 0 & y-1 & -y & 1   \\ \hline
 \end{array}$

 \end{multicols}
 
The marked point is to the right of the $5^{th}$ crossing, following the dotted path we see that we can see that $w_5=y^{-2}$.  Checking that  $\mu (1) =1$, $\mu(2)= 4$, $rot(1)=1$, $rot(2)=3$, and $M_5^5=x(y-1)(1-y+y^2)$  gives the MVA of the link as $xy(1-y+y^2)$.

\section{The Multivariable Alexander Polynomial as a Finite Type Invariant}

In \cite{Mur} Murakami uses the following to show that under appropriate change of variables the MVA is of finite type.

A knot with a singular point is really just the difference between the two following knots.  Note that the added kink is just to make the number and labeling of arcs consistent, and does not change the knot.  The  two matrices constructed as in Section 1 will differ in only the two rows indexed by $a_1$ and $a_2$.  Note that the kink doesn't change the normalization as it adds one to $\mu(k)$ and $rot(k)$.  The only difference in the normalization coefficient for these two knots is the number of over crossings that adds an extra factor of $t_i^{\frac{1}{2}}$.  By multiplying  a matrix row by that change we can see that the matrices differ only in the following rows:

\begin{center}
\includegraphics[width=80mm]{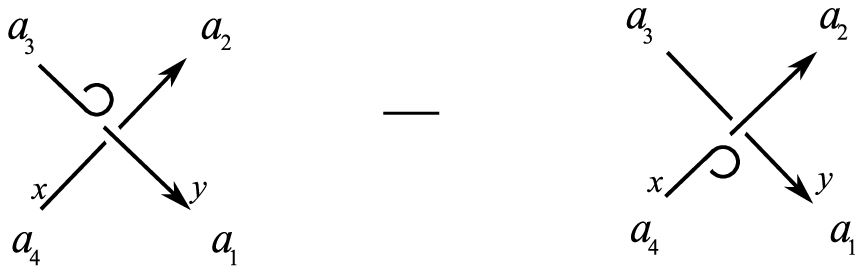}
\end{center}
\begin{multicols}{2}
$\begin{array}{c|c c c c }
 & a_1 & a_2 & a_3 & a_4\\  \hline
 
 a_1 & y^
{-{\frac{1}{2}} }& 0 & -y^{-{\frac{1}{2}} } & 0 \\
 a_2 & 1-x & y & 0 & -1 \\
\end{array}$

$\begin{array}{c|c c c c }
 & a_1 & a_2 & a_3 & a_4\\  \hline
 
 a_1 & 1 & y-1 & -x & 0 \\
 a_2 & 0 & x^{-{\frac{1}{2}} } & 0 & -x^{-{\frac{1}{2}} } \\
\end{array}$

\end{multicols}

 Using row operations we see that the resulting difference will contain the following rows:

  {\center \noindent $\begin{array}{c|c c c c }
 & a_1 & a_2 & a_3 & a_4\\  \hline
 
 a_1^{\prime} & y^{-{\frac{1}{2}} }  & x^{-{\frac{1}{2}} }  & -y^{-{\frac{1}{2}} }  & -x^{-{\frac{1}{2}} }  \\
 a_2 ^{\prime}& 1 -x^{{\frac{1}{2}} }& y-y^{{\frac{1}{2}} }  & -x+x^{{\frac{1}{2}} } & -1+y^{{\frac{1}{2}} }\\
\end{array}.$ }

\

Now if we use the substitution  $e^{t_1}=x$, $e^{t_2}=y$  and expand  $e^{t_i}$ as a power series, omitting terms of degree greater than 1 we get;

 {\center \noindent $\begin{array}{c|c c c c }
 & a_1 & a_2 & a_3 & a_4\\  \hline
 
 a_1^{\prime} & 1-\frac{1}{2}t_2-\ldots  & 1-\frac{1}{2}t_1- \ldots  & -1 +\frac{1}{2}t_2+\ldots & -1+\frac{1}{2}t_1+\ldots \\
 a_2 ^{\prime}& -\frac{1}{2}t_1-\ldots& \frac{1}{2}t_2+\ldots& -\frac{1}{2}t_1-\ldots&\frac{1}{2}t_2+\ldots\\
\end{array}$ }

\

Then for each double point we get a row that has all terms of degree equal to or higher than 1 so  the determinant of the matrix has  degree of at least m.  When we apply the same change of variables to the normalization coefficient, we notice that dividing by $x_i-1$ will lower the index of $t_i$ by one.  Since the rest of the coefficient will have leading term 1,  the resulting polynomial will have degree at least  $m-1$.  Hence the term of degree $m$ is an $m+1$ type invariant.

\section{The Weight System}

We define a function from  coloured chord diagrams to $\mathbb{Z}[t_1, ... , t_n]$.  Given a coloured chord diagram, choose a marked point on one of the arcs.  This divides that arc into two; we then label the arcs.    The chords are indexed by the two exiting arcs.  We now construct a matrix $M(D)$ whose rows and columns are indexed by the arcs, as follows. (Notice if you change roles of the two sides, you swap the two rows below and multiply by $-1$ which leaves the determinant unchanged).

\begin{multicols}{2}

\begin{center}
\includegraphics[width=40mm]{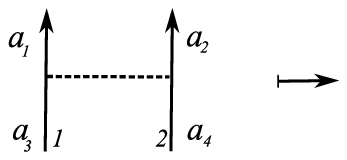}
\end{center}

\ \ 
\begin{tabular}{|l | l l l l |} \hline

& $a_1$ & $a_2$ & $a_3$ & $a_4$  \\    \hline

$a_1$ & $\frac 1 2 $  &  $ \frac{1}{2} $ & $\frac{-1}{2}$  &  $ -\frac{1}{2} $ \\
$a_2$ & $-t_2$  &  $ t_1 $ & $-t_2$  &  $ t_1 $ \\ \hline
\end{tabular}

\end{multicols}

Note that the row corresponding to the marked point is empty. As before $M_i^i(D)$ is the matrix $M(D)$ with the $i^{th}$ row and column removed.

\begin{definition}{The MVA weight is}

 $$\Delta (D) =\frac{ \det(M_i^i(D))}{t_i}$$
 
 where $i$ is the colour of the marked point.
 
 \end{definition}

To aid in later sections we can extend the definition to allow us to divide arcs at non chords.  To divide the arc $a$ into $a_1$ and $a_2$ the the following row would need to be added to the matrix.  Note that this will not change the value of the determinant.

\begin{multicols}{2}

\begin{center}
\includegraphics[width=27mm]{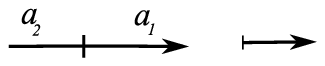}
\end{center}

$\begin{array}{|c|c c |}\hline
  &  a_1  &  a_2  \\
  \hline   a_1 &  1  & -1  \\ \hline

\end{array}$
\end{multicols}

\begin{theorem}
The above function defines a weight system for the Multivariable Alexander Polynomial under the substitution $x_i=e^{t_i}=1+t_i+\dots$
\end{theorem}

A proof will follow a sample computation.

\section{A Worked Example}

\includegraphics[width=65mm]{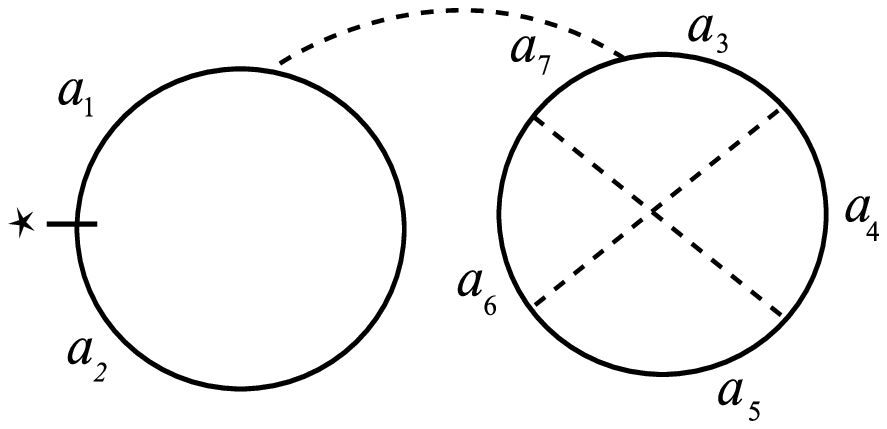}

On the above diagram we have chosen a marked point $\star$ and labeled the arcs. Note that the left circle is coloured by $x$ and the right by $y$.  We now construct a matrix indexed by those arcs.  Here is one of the chords and the corresponding matrix entries.

\begin{multicols}{2}

\begin{center}
\includegraphics[width=40mm]{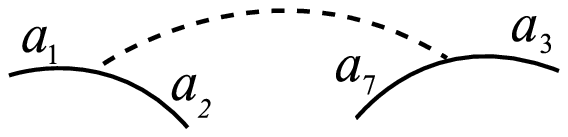}
\end{center}

\ \ 
\begin{tabular}{|l | l l l l |} \hline

& $a_1$ & $a_2$ & $a_3$ & $a_7$  \\    \hline

$a_2$ & $\frac {-1} {2} $  &  $ \frac{1}{2} $ & $\frac{1}{2}$  &  $ \frac{-1}{2} $ \\
$a_3$ & $-t_2$  &  $ -t_2 $ & $t_1$  &  $ t_1 $ \\ \hline
\end{tabular}

\end{multicols}

We can use the rest of the chords to fill in the remainder of the table as follows.

\begin{equation*}
\begin{array}{|c | c c c c c c c  | } \hline

  & a_{1} &  a_{2} &  a_{3} & a_{4} &  a_{5} &  a_{6} & a_{7} \\  \hline
a_{1} &  & & & & & &\\
 a_{2} &  -\frac{1}{2} & \frac{1}{2} & \frac{1}{2} &  &  &  &  -\frac{1}{2} \\

 a_{3} & -t_2 & -t_2 & t_1 &  &  &  &  t_1  \\
a_{4}  &  &  &  t_2 &  t_2 &  -t_2 &  -t_2 &    \\
a_{5} & & &  &  -t_2  &  -t_2 & t_2 &  t_2 \\
  a_{6}  &  &  &  -\frac{1}{2}  &  \frac{1}{2}  &  -\frac{1}{2}  & \frac{1}{2} &   \\
 a_{7} &  &  &  &  -\frac{1}{2}  & \frac{1}{2} & -\frac{1}{2} & \frac{1}{2} \\  \hline
\end{array}
\end{equation*}

To calculate the weight of this diagram, we simply remove the marked row and column ($a_{1}$), take the determinant, divide by $t_1$ to get  $-t_2^2$.

\section{Proof of Theorem}

Suppose we have an $m$ singular link and its underlying chord diagram.  We can label the arc of the link as shown,  with $x_i$ indexing the arcs on the underlying chord diagram and $x_{i,j}$ the arcs on the link.  Then we mark an arc to correspond to the row and column that we will delete.  The crossings are indexed by the under arc leaving it.

\begin{center}
\includegraphics[width=120mm]{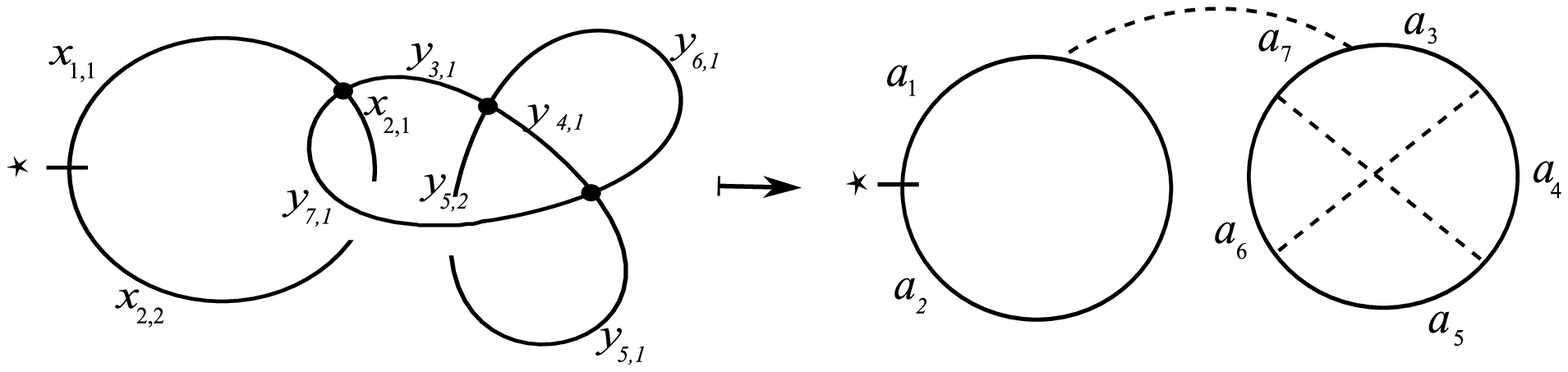}
\end{center}

As in the proof of finite type, we will be using the substitution $x_i=e^{t_i}=1+t_i+\ldots$.  Since we are interested in the terms of order $m-1$, we take a row of order $1$ in the rows corresponding to the double points, and terms of order zero elsewhere. We note that up to order $0$ the rows corresponding to over and undercrossings are the same after this substitution.  So to we will create  a matrix with the following rows;

\begin{multicols}{2}

\begin{center}
\includegraphics[width=60mm]{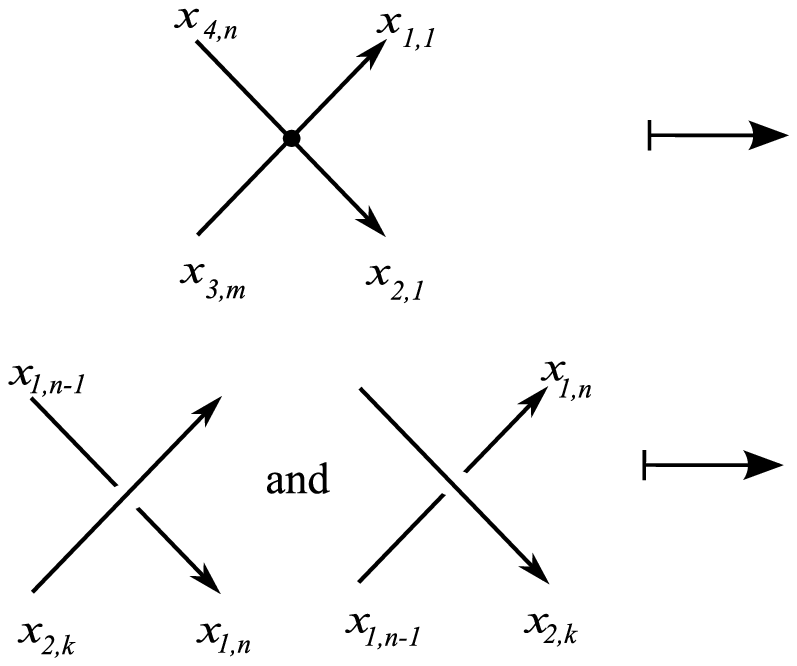}
\end{center}

\

\vspace{5mm}
\begin{tabular}{|l | l l l l |} \hline

& $x_{1,1}$ & $x_{2,1}$ & $x_{3,n}$ & $x_{4,m}$  \\    \hline

$x_{1,1}$ & $\frac {1} {2} $  &  $ \frac{1}{2} $ & $\frac{-1}{2}$  &  $ \frac{-1}{2} $ \\
$x_{2,1}$ & $t_2$  &  $ -t_1 $ & $t_2$  &  $ -t_1 $ \\ \hline
\end{tabular}

\vspace{13mm}

\begin{tabular}{|l | l l l |} \hline

& $x_{1,n}$ & $x_{1,n-1}$ & $x_{2,k}$   \\    \hline

$x_{1,n}$ & $1 $  &  $-1 $& 0   \\
 \hline
\end{tabular}

\end{multicols}

 This gives us a matrix of mostly $1$ and $-1$'s.  The only variables occur in the rows $x_{i,1}$.  The matrix will be almost block diagonal.  The only exception is the top row of each  block.  Each block will also have a specific form, as shown below, with zeros omitted.\\

\begin{tabular}{l|c c c c c c  c | }

   & $x_{1 ,1}$ & $x_{1,2} $ & $ x_{1,3} $ & $ ...$ & $x_{1,n-1}$ & $x_{1,n}$ \\  \hline
   
 $x_{1 ,1}$ &   *  &   &  &   &    &  *  \\
 
  $x_{1 ,2}$ &   -1  &  1  &     &     &    &    \\
  
   $x_{1 ,3}$ &     &    -1 &   1  &     &    &    \\
   
    \vdots &     &     &   $\ddots$ &  $\ddots$   &    &    \\
    
   $x_{1 ,n-1}$ &     &     &     &   -1  &  1  &    \\
   
    $x_{1 ,n}$ &     &     &     &     &  -1  &   1 \\

\end{tabular}

Then for all of the blocks except for the marked one,  we can switch to the following using column operations.\\

\begin{tabular}{l|c c c c c c  c | }

   & $x_{1 ,1}$ & $x_{1,2} $ & $ x_{1,3} $ & $ ...$ & $x_{1,n-1}$ & $x_{1,n}$ \\  \hline
   
 $x_{1 ,1}$ & $\Sigma  * $ &     &     &     &    &  *  \\
 
  $x_{1 ,2}$ &   0  &  1  &     &     &    &    \\
  
   $x_{1 ,3}$ &   0  &   0&   1  &     &    &    \\
   
    \vdots &     &     &   $\ddots$ &  $\ddots$   &    &    \\
    
   $x_{1 ,n-1}$ &     &     &     &  0 &  1  &    \\
   
    $x_{1 ,n}$ &    &     &     &     & 0  &   1 \\

\end{tabular}

Note that rows $x_{1 ,2}$  to $x_{1 ,n}$ are zero except with a $1$ on the diagonal, so we may delete those rows and columns without changing the determinant.  But we must remember that the $x_{1 ,1}^{th}$ column is replaced with the sum of columns $x_{1 ,1}$  and $x_{1 ,n}$.  That sum will be the contribution of the arcs $x_{1,1}$ and $ x_{1,n}$  from double points, note that those arcs map to $x_1$ in the chord diagram.

In the case that one of the arcs is  marked, we must remember that we will be deleting that column before we take the determinant.  We will always mark an arc that is not near a double point.  Let us look at the block where one of the arcs is marked; with out loss of generality we choose $x_{1,2} $.  Then when we reduce the matrix we get:

\begin{tabular}{l|c c c c   c | }

   & $x_{1 ,1}$ & $ x_{1,3} $ & $ ...$ & $x_{1,n-1}$ & $x_{1,n}$ \\  \hline
   
      $x_{1 ,1}$ & *    &  &   &    &  *  \\
      $x_{1 ,3}$ &    0 &   1  &     &    &    \\
   
    \vdots &       &   $\ddots$ &  $\ddots$   &    &    \\
    
   $x_{1 ,n-1}$ &         &     &   -1  &  1  &    \\
   
    $x_{1 ,n}$ &         &     &     &  -1  &   1 \\

\end{tabular}

Which can then be changed using row and column operations to:

\begin{tabular}{l|c c c c  c | }

   & $x_{1 ,1}$ & $ x_{1,3} $ & $ ...$ & $x_{1,n-1}$ & $x_{1,n}$ \\  \hline
   
 $x_{1 ,1}$ &   *  &   &   &    &  *  \\

   $x_{1 ,3}$ &    0 &   1  &     &    &    \\
   
    \vdots &       &   $\ddots$ &  $\ddots$   &    &    \\
    
   $x_{1 ,n-1}$ &         &     &  0  &  1  &    \\
   
    $x_{1 ,n}$ &         &     &     &  0  &   1 \\

\end{tabular}

Then once again we can delete the rows and columns corresponding to $ x_{1,3} $ , $ ...$ , $x_{1,n-1}$.  
But note that there is no change in the $x_{1,1}^{th}$ column.  So the contribution of the  $x_{1,n}^{th}$ column is lost, this is equivalent to deleting the column associated to the marked chord.

 This is the matrix that was described in Theorem 1. \\

\section{Testing Known Relations}

As a test of this method, we wish to show that this function obeys certain known relations, such as the 4T relation and the recursive relations in Murakami's paper.  Each of these relations involve chord diagrams that differ only at one small area.  This means that the matrix formed in computing the invariant will differ in only  one block.  The following lemma shows that the relations need only hold among certain minors of that block.

\begin{lemma}[A linear algebra lemma]
 Let B be an $n\times n$ block matrix of the following form, where A is an $m \times l$ matrix, and  $k<l$.

\begin{center} 
\includegraphics[width=30mm]{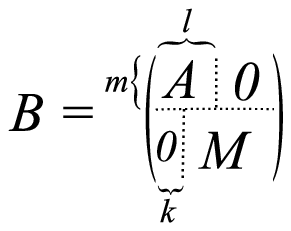}
\end{center}

\mbox{Then} 
\begin{equation*} \det(B) = \sum _{k<i_1<...<i_{n-k}} (\pm ) \det (A^{1,2,\ldots,k,i_1,...,i_{n-k}})\det (M_{i_1,...,i_{n-k}}),
\end{equation*}

where $A^{1,2,\ldots,k,i_1,...,i_{n-k}}$ refers to the matrix formed from the $1,2,\ldots,k,i_1,...,i_{n-k}$ columns of A and $M_{i_1,...,i_{n-k}}$ refers to M with the $i_1,...,i_{n-k}$ columns removed.

\end{lemma}

 \noindent{\bf Proof:}

Expand B along the rows of A using cofactor expansion, and collect like terms.\\

\begin{lemma}
Let $B_i$ be matrices as in the previous lemma, with the same $M$ but different $A_i$.  To show a relation of the form

$$\sum a_i\det (B_i) =0,$$

it is enough to show the relation holds on the appropriate minors the $A_i$'s.

\end{lemma}

 \noindent{\bf Proof:}
 
 Follows from previous lemma.\\

When we wish to show relations between chord diagrams differing in a region, we note that the matrices used to compute their weight will be of the above form.  Using this lemma let us show the following relation holds.

\includegraphics[width=50mm]{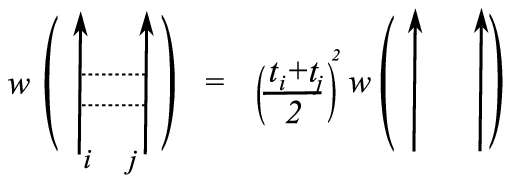}

To use the previous lemma we make both sides of the relation have the same number of arcs, so we divide the arcs in the left term in three.  The matrices involved will be as follows.

\begin{multicols}{2}

$
\begin{array}{|c c c c c c|c | } \hline

  \frac{1}{2} & \frac{1}{2} & 0 & 0 & -\frac{1}{2} & -\frac{1}{2}  & 0 \\ 
   
 -t_2 & t_1 & 0 & 0 & -t_2 & -t_1  & 0 \\  
  -\frac{1}{2} & -\frac{1}{2} & \frac{1}{2} & \frac{1}{2} & 0 & 0  & 0 \\  
  
   -t_2 & t_1 & -t_2 & t_1 & 0 & 0 & 0  \\  \hline
   
    0 & 0& a_1 & a_2 & a_3 & a_4  & M \\  \hline

\end{array}
$

$
\begin{array}{|c c c c c c | c| } \hline

  1 & 0 & 0 & 0 & -1 & 0 & 0 \\ 
   
 0 & 1 & 0 & 0 & 0 & -1 & 0 \\  
 
 -1 & 0 & +1 & 0 & 0 & 0 & 0 \\ 
  
   0 & -1 & 0 & 1 & 0 & 0 & 0 \\\hline
 0 & 0& a_1 & a_2 & a_3 & a_4 & M \\  \hline

\end{array}$

\end{multicols}

For the first matrix in the example above the minors in question are as follows;

 \begin{multicols}{2}

$$\begin{array}{|c c c c| }
 \frac 1  2 & \frac 1  2 & 0 & 0 \\ 
   -t_2 & t_1 & 0 & 0  \\  
  -\frac 1  2 & -\frac 1  2 & \frac 1  2 & \frac 1  2  \\  
   -t_2 & t_1 & -t_2 &  t_1 \\
\end{array} = \left(\frac{t_1 + t_2}{2}\right)^2$$

$$\begin{array}{|c c c c |} 

  \frac 1  2 & \frac 1  2 & 0  & -\frac 1  2   \\ 
   
 -t_2 & t_1 & 0  & -t_2  \\  
  -\frac 1  2 & -\frac 1  2 & \frac 1  2  & 0    \\  
  
   -t_2 & t_1 & -t_2 & 0   \\

\end{array} = 0$$

$$\begin{array}{|c c c c   | } 

  \frac 1  2 & \frac 1  2 & 0 & -\frac 1  2   \\ 
   
 -t_2 & t_1 & 0 & -t_1   \\  
  -\frac 1  2 & -\frac 1  2 & \frac 1  2 & 0   \\  
  
   -t_2 & t_1 & -t_2 & 0  \\

\end{array} = -\left(\frac{t_1+t_2}{2}\right)^2$$

$$\begin{array}{|c c c c | }

  \frac 1  2 & \frac 1  2 &  0 & -\frac 1  2  \\ 
   
 -t_2 & t_1 & 0 & -t_2 \\  
  -\frac 1  2 & -\frac 1  2  & \frac 1  2 & 0  \\  
  
   -t_2 & t_1  & t_1 & 0 \\

\end{array} = \left(\frac{t_1+t_2}{2}\right)^2$$

$$\begin{array}{|c c c c | }

  \frac 1  2 & \frac 1  2& 0  & -\frac 1  2 \\ 
   
 -t_2 & t_1& 0  & -t_1  \\  
  -\frac 1  2 & -\frac 1  2 & \frac 1  2  & 0  \\  
  
   -t_2 & t_1 & t_1 & 0  \\

\end{array} = 0$$

$$\begin{array}{|c c c c | } 

  \frac 1  2 & \frac 1  2 & -\frac 1  2 & -\frac 1  2 \\ 
   
 -t_2 & t_1 & -t_2 & -t_1   \\  
  -\frac 1  2 & -\frac 1  2 & 0 & 0   \\  
  
   -t_2 & t_1 &  0 & 0   \\  

\end{array} = \left( \frac{t_1+t_2}{2} \right) ^2.$$

\end{multicols}

The corresponding minors for the second matrix are, {$1, 0, -1, 1, 0, 1$}.  You can see that the first are  precisely $\left(\frac{t_1+t_2}{2}\right)^2$ times those of the second,  so the relation  holds.

If we cross the chords we get the following two relations.\\

\includegraphics[width=100mm]{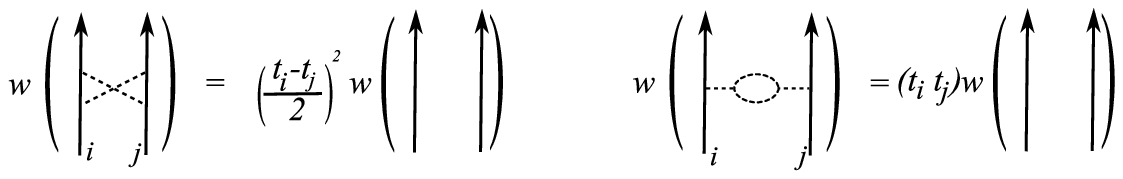} \\

This allows us to derive the following  useful relations.  \\

\includegraphics[width=100mm]{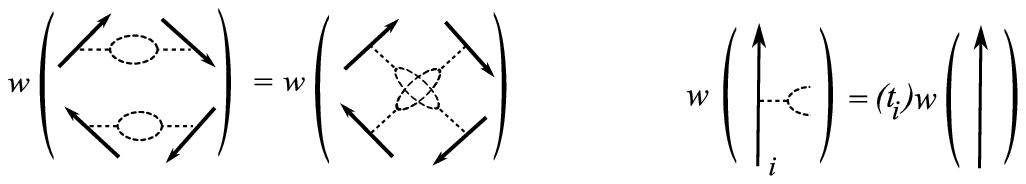}

The first relation is a direct consequence of the previous.  Since it doesn't matter how these components connect we can omit the connections from our diagrams.  We can draw elements as in the second relation with the understanding that there will always be an even number of such elements in a diagram.

\section{Testing Relations Using Mathematica}

Obviously for more complicated relations we are unable to compute by hand all of the minors involved.  The solution to this is to allow Mathematica to do all of the computations for us.  For example here is the computation that shows the 4Y relation holds.\\

\includegraphics[width=110mm]{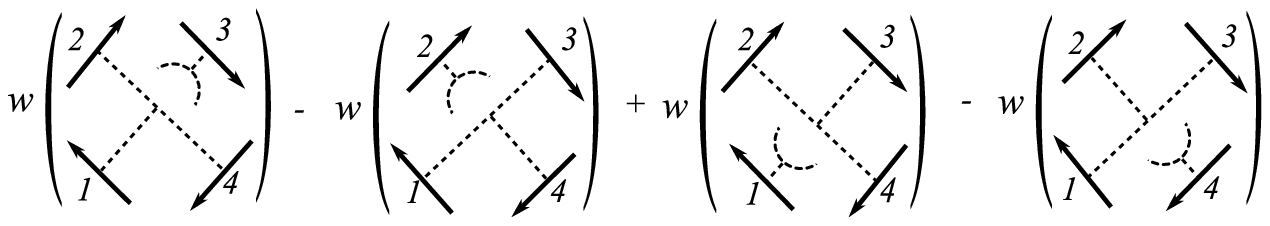}

Using the previous relation we can rewrite the relation as.\\

\includegraphics[width=110mm]{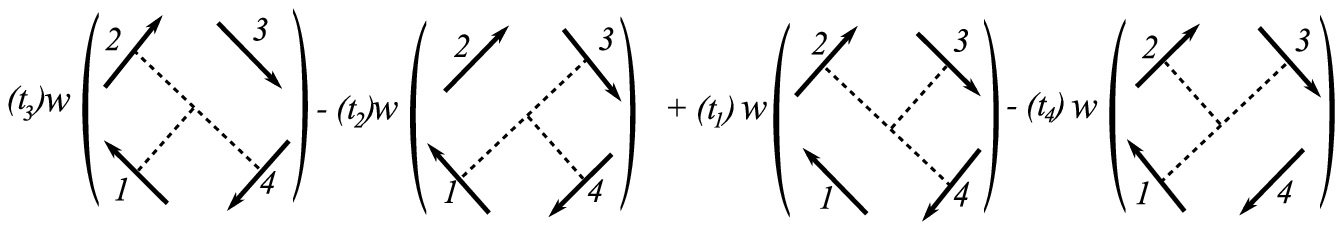}

In the following calculation, we have used the substitution $x=t_1$, $y=t_2$, $z=t_3$ and $w=t_4$, and  a common factor of $\frac{1}{8}$ has been factored from all of the matrices.

\begin{multicols}{2}

\begin{verbatim}
Mat[1]={
      {1,1,0,0,0,-1,-1,0,0},
      {-x,y,0,0,0,y,-x,0,0},
      {-1,0,1,0,1,0,0,0,-1},
      {0,0,0,-1,0,0,0,1,0},
      {-w,0,-w,0,y,0,0,0,y}};

Mat[2]={
      {w,0,0,0,-y,0,w,0,-y},
      {-1,1,1,0,0,-1,0,0,0},
      {x,-y,x,0,0,-y,0,0,0},
      {0,0,0,-1,0,0,0,1,0},
      {1,0,0,0,1,0,-1,0,-1}};

Mat[3]={
      {z,0,0,-w,0,0,0,-w,z},
      {-1,1,0,0,1,-1,0,0,0},
      {0,0,-1,0,0,0,1,0,0},
      {1,0,0,1,0,0,0,-1,-1},
      {x,-w,0,0,x,-w,0,0,0}};

Mat[4]={
      {x,-w,0,0,0,-w,0,0,x},
      {1,1,0,0,0,-1,0,0,-1},
      {0,0,-1,0,0,0,1,0,0},
      {-1,0,0,1,1,0,0,-1,0},
      {z,0,0,-w,z,0,0,-w,0}};

Mat[5]={
      {z,0,0,-w,0,0,0,-w,z},
      {0,-1,0,0,0,1,0,0,0},
      {-y,0,w,0,-y,0,w,0,0},
      {1,0,0,1,0,0,0,-1,-1},
      {-1,0,1,0,1,0,-1,0,0}};

Mat[6]={
      {1,0,1,0,0,0,-1,0,-1},
      {0,-1,0,0,0,1,0,0,0},
      {-y,0,w,0,0,0,w,0,-y},
      {-1,0,0,1,1,0,0,-1,0},
      {z,0,0,-w,z,0,0,-w,0}};

Mat[7]={
      {1,1,0,0,0,-1,-1,0,0},
      {-x,y,0,0,0,y,-x,0,0},
      {z,0,z,-y,0,0,0,-y,0},
      {-1,0,1,1,0,0,0,-1,0},
      {0,0,0,0,-1,0,0,0,1}};
      
Mat[8]={
      {z,0,0,-y,0,0,z,-y,0},
      {-1,1,1,0,0,-1,0,0,0},
      {x,-y,x,0,0,-y,0,0,0},
      {1,0,0,1,0,0,-1,-1,0},
      {0,0,0,0,-1,0,0,0,1}};
\end{verbatim}
\end{multicols}
\begin{verbatim}
z*(Minors[Mat[1],5]-Minors[Mat[2],5])-y*(Minors[Mat[3],5]-
Minors[Mat[4],5])+x*(Minors[Mat[5],5]-Minors[Mat[6],5])-
w*(Minors[Mat[7],5]-Minors[Mat[8],5]);

L=Simplify[%];

Take[Transpose[L],70]

{{0},{0},{0},{0},{0},{0},{0},{0},{0},{0},{0},{0},{0},{0},{0},{0},
{0},{0},{0},{0},{0},{0},{0},{0},{0},{0},{0},{0},{0},{0},{0},{0},
{0},{0},{0},{0},{0},{0},{0},{0},{0},{0},{0},{0},{0},{0},{0},{0},
{0},{0},{0},{0},{0},{0},{0},{0},{0},{0},{0},{0},{0},{0},{0},{0},
{0},{0},{0},{0},{0},{0}}
\end{verbatim}

From our lemma, since the relation holds for the minors, it holds for the determinants of the original matrices.  So the relation holds for the weight system.

\section{Additional Relations}

In \cite{gl} additional relations where shown to hold for the weight system of the Alexander polynomial, these can be shown to hold for the MVA weight system as well.

\

No internal vertices

\

\includegraphics[width=60mm]{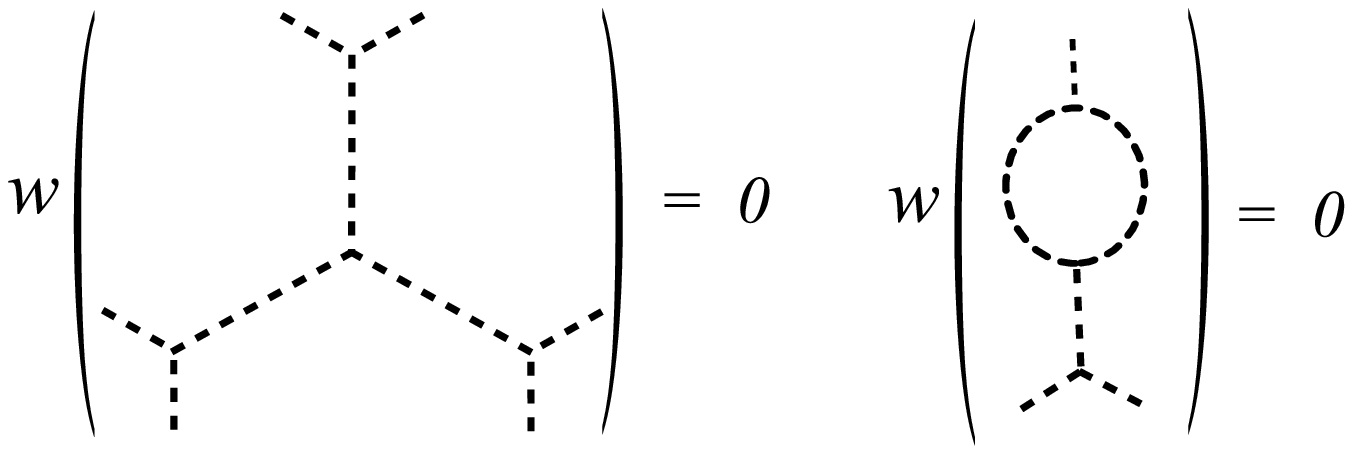}

\

IHX

\

\includegraphics[width=80mm]{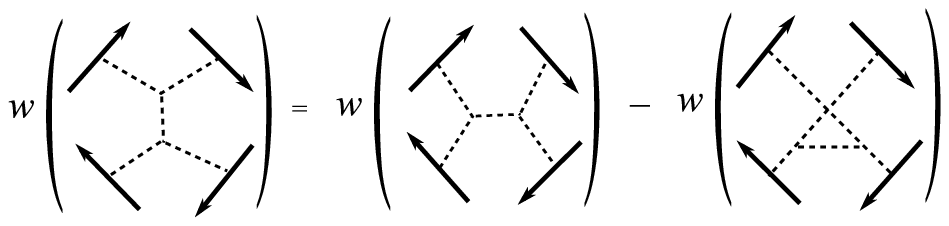}

\end{document}